\documentclass[a4paper,10pt,oneside]{article}
\usepackage{amssymb}
\usepackage{amsmath}
\usepackage{amsthm}
\title{\textbf{Fractional dynamical systems defined on fractional jet bundle and applications in economics }}
\author{Mihai Boleantu \\  Dept. of Economics,
West University of Timisoara \\
email: mihai.boleantu@fse.uvt.ro and \\ mihaiboleantu01@yahoo.com}

\begin{document}
\date{} %
\maketitle
\begin{abstract}

Using Caputo fractional  derivative of order $\alpha $ we build the fractional jet bundle of order $\alpha $ and its main geometrical structures. Defined on that bundle, some fractional dynamical systems with applications to economics are studied.

\end{abstract}
\newtheorem{theorem}{Theorem}
\newtheorem{proposition}[theorem]{Proposition}
\renewcommand{\arraystretch}{1.5}

\section{Introduction}

Fractional integration and fractional differentiation are generalizations of the notions of integer-order integration and differentiation, and include \textit{n}-th derivatives and \textit{n}-fold integrals as a particular case (\textit{n}  denotes an integer number). Fractional calculus has been used successfully in various fields of science and engineering. These applications include classical and quantum mechanics, field theory, optimal control and economics. In the past few decades many authors \cite{6}, \cite{8}, \cite{10} have pointed out that fractional-order models are more appropriate than integer-order models for the description of memory and hereditary properties of various  processes. This is the main advantage of fractional derivatives in comparison with classical integer-order models in which such effects are, in fact, neglected. For example, the advantages of fractional derivatives become apparent in modeling mechanical and electrical properties of real materials, as well as in the description of rheological properties of rocks (see \cite{2}, \cite{6}). Also, in economics, fractional calculus has been used for deriving an analytical model of the tax version of the Fisher equation that incorporates a memory function for stock prices and inflation rates (see \cite{4}). \\
\indent The problems have been formulated mostly in terms of two types of fractional derivatives, namely Riemann-Liouville (RL) and Caputo \cite{2}, \cite{9}. Among mathematicians, RL fractional derivatives have been popular largely because they are amenable to many mathematical manipulations. However, the RL derivative of a constant is not zero, and in many applications it requires fractional initial conditions which are generally not specified. In contrast, Caputo derivative of a constant is zero, and a fractional differential equation defined in terms of Caputo derivatives require standard boundary conditions. For these reasons, Caputo fractional derivatives have been popular among engineers and scientists. \\ \indent A simple and really geometric interpretation of several types of fractional-order integration is given in \cite{9}. Based on this, a physical interpretation of the Riemann-Liouville fractional integration is proposed in terms of inhomogeneous and changing (nonstatic, dynamic) time scale, and  it is shown that the Caputo fractional derivative has the same physical interpretation as the Riemann- Liouville fractional derivative. \\
\indent In this paper, using Caputo fractional derivative, the fractional jet fibre bundle is built  on a differentiable manifold and its main geometric structures are emphasized. Some significant examples from economics are presented. In section 2, some useful properties of the Caputo fractional derivative are recalled and the fractional osculator bundle of order \textit{k} is described. In section 3, the fractional jet fibre bundle is defined and the fractional Euler-Lagrange equations together with the fractional Hamilton equations are established. In section 4 two fractional economic models are studied.
\section{The fractional osculator bundle of order \textit{k}
on a differentiable manifold}
\subsection{Caputo fractional derivative of  order $\alpha $} Let $x:[a,b]\to {\mathbb{R}}$ be a smooth function, $\alpha \in {\mathbb{R}}$, $\alpha \ge 0$, and $m\in {\mathbb N}^{*} $ such that $m-1<\alpha <m$. The left (right) \textit{Caputo} \textit{fractional derivative} of $x$ \cite{2} is the function

\begin{equation} \label{ZEqnNum915072} \begin{array}{l} {{}_{a} D_{t}^{\alpha } x(t)=\frac{1}{\Gamma (m-\alpha )} \int _{a}^{t}\frac{x^{(m)} (s)}{(t-s)^{\alpha +1-m} } ds } \\ {{}_{t} D_{b}^{\alpha } x(t)=\frac{1}{\Gamma (m-\alpha )} \int _{t}^{b}\frac{(-1)^{m} x^{(m)} (s)}{(t-s)^{\alpha +1-m} } ds ,} \end{array} \end{equation}

\begin{flushleft}
where $x^{(m)} (s)=\frac{d^{m} }{ds^{m} } x(s)$ and $\Gamma $ is the gamma function of Euler. The following properties result from \eqref{ZEqnNum915072} (see \cite{3}).

1. ${}_{a} D_{t}^{\alpha } (c_{1} x_{1} (t)+c_{2} x_{2} (t))=c_{1} {}_{a} D_{t}^{\alpha } x_{1} (t)+c_{2} {}_{a} D_{t}^{\alpha } x_{2} (t),$ $\forall c_{1} ,c_{2} \in {\mathbb{R}}$, $\forall x_{1} ,x_{2} :[a,b]\to {\mathbb{R}}$.

2. $D_{t}^{1} {}_{a} D_{t}^{\alpha } x(t)={}_{a} D_{t}^{\alpha +1} x(t)+\frac{t^{m-\alpha -1} }{\Gamma (m-\alpha )} x^{(m)} (a)$, $D_{t}^{1} =\frac{d}{dt} $.

3. If $\left\{\alpha _{n} \right\}_{n\ge 0} $ is a real number sequence with $\mathop{\lim }\limits_{n\to \infty } \alpha _{n} =p\in {\mathbb N}^{*} $, then $\mathop{\lim }\limits_{n\to \infty } \left({}_{a} D_{t}^{\alpha _{n} } x\right)(t)=D_{t}^{1} x(t)$.

4. (a) If $x(t)=c$, $\forall t\in [a,b]$, $c\in {\mathbb R}$, then ${}_{a} D_{t}^{\alpha } x(t)=0$.

    (b) If $x(t)=t^{\gamma } $, $\forall t\in [a,b]$, ${}_{a} D_{t}^{\alpha } x(t)=\frac{t^{\gamma -\alpha } \Gamma (1+\gamma )}{\Gamma (1+\gamma -\alpha )} $.

5. If $x_{1} ,x_{2} :[a,b]\to {\mathbb{R}}$ are analytic functions, then

\[{}_{a} D_{t}^{\alpha } (x_{1} x_{2} )(t)=\sum _{k=0}^{\infty }\binom{\alpha}{k} {}_{a} D_{t}^{\alpha -k} x_{1} (t)D_{t}^{k} x_{2} (t),\]

where $D_{t}^{k} =\frac{d}{dt} \circ ...\circ \frac{d}{dt} $.
\end{flushleft}
\begin{flushleft}
6.  $\int _{a}^{b}x_{1} (t)\left({}_{a} D_{t}^{\alpha } x_{2} (t)\right)dt =-\int _{a}^{b}x_{2} (t)\left({}_{a} D_{t}^{\alpha } x_{1} (t)\right)dt $.

7. If $x:[a,b]\to {\mathbb R}$ is an analytic function and $0\in (a,b)$ then

\begin{equation} \label{2)} x(t)=\sum _{a=0}^{\infty }\frac{t^{\alpha a} }{\Gamma (1+\alpha a)}  {}_{a} D_{t}^{\alpha a} x(t)\left|_{t=0} \right. . \end{equation}
\end{flushleft}

\subsection{Higher-order fractional osculator bundle }

 Let $\alpha \in (0,1)$ be fixed and $M$ a differentiable manifold of dimension $n$. Two curves $c_{1} ,c_{2} :I\to {\mathbb R}$ with $c_{1}(0)=c_{2}(0)=x_{0} \in M$, $0\in I$, have a fractional contact $\alpha $ of order $k\in {\mathbb N}^{*} $ in $x_{0} $, if for any $f\in {\mathcal{F}}(U)$, $x_{0} \in U$, $U$ a chart on $M$, it holds that

\begin{equation} \label{ZEqnNum982126} D_{t}^{\alpha a} (f\circ c_{1} )\left|_{t=0} \right. =D_{t}^{\alpha a} (f\circ c_{2} )\left|_{t=0} \right. , \end{equation}
where $a=\overline{1,k}$ and $D_{t}^{\alpha a} ={}_{0} D_{t}^{\alpha a} $.

 The set of equivalence classes defined by \eqref{ZEqnNum982126} is called the \textit{k}-\textit{osculator fractional space} at $M$ in $x_{0} $ and it will be denoted by $Osc_{x_{0} }^{\alpha k} (M)$. If the curve $c:I\to M$ is given by $x^{i} =x^{i} (t)$, $t\in I$, $i=\overline{1,n}$, in the chart $U$, then the class $[c]_{x_{0} }^{\alpha k} \in Osc_{x_{0} }^{\alpha k} (M)$ is given by

\begin{equation} \label{4)} \begin{array}{l} {x^{i} (t)=x^{i} +\sum _{a=1}^{k}t^{\alpha a} y^{i(\alpha a)}  ,\; \; t\in (-\varepsilon ,\varepsilon ),} \\ {y^{i(\alpha a)} =\frac{1}{\Gamma (1+\alpha a)} D_{t}^{\alpha a} x^{i} (t)\left|_{t=0} \right. ,} \end{array} \end{equation}

where $i=\overline{1,n}$, $a=\overline{1,k}$ and $x^{i} =x{}^{i}(0)$. The \textit{fractional osculator bundle of order k} is the bundle $(Osc^{\alpha k} (M),\; M)$ where $Osc^{\alpha k} (M)=\bigcup _{x\in M}Osc_{x}^{\alpha k} (M) $ and $\pi ^{\alpha k} :Osc^{\alpha k} (M)\to M$ is defined by $\pi ^{\alpha k} ([c]_{x}^{\alpha k} )=x$, $\forall [c]_{x}^{\alpha k} \in Osc^{\alpha k} (M)$.

 For $f\in {\mathcal{F}}(U)$, the fractional derivative of order $\alpha $, $\alpha \in (0,1)$, with respect to the variable $x^{i} $ is defined by

\[D_{x^{i} }^{\alpha } f(x)=\frac{1}{\Gamma (1-\alpha )} \int _{a^{i} }^{x^{i} }\frac{\partial f(x^{1} ,...,x^{i-1} ,s,x^{i+1} ,...,x^{n} )}{\partial x^{i} } \frac{1}{(x^{i} -s)^{\alpha } } ds, \]

where $x\in U_{ab} =\{ x\in U,\; a^{i} \le x^{i} \le b^{i} ,\; i=\overline{1,n}\} $, $x^{i} $, $i=\overline{1,n}$, are the coordinate functions on $U$ and $\left\{\frac{\partial }{\partial x^{i} } \right\}$, $i=\overline{1,n}$, is the canonical base of the vector fields on $U$.

 Using the fractional exterior differential \cite{1} $d^{\alpha } :{\mathcal{F}}(U)\to {\mathcal{D}}^{1} (U)$ given by

\begin{equation} \label{5)} d^{\alpha } =d(x^{j} )^{\alpha } D_{x^{j} }^{\alpha } , \end{equation}

where $(x^{j} )^{\alpha } \in {\mathcal{F}}(U)$ and ${\mathcal{D}}^{1} (U)$ is the module of the differential 1-forms on $U$ with the canonical base $\{ dx^{i} \} $, $i=\overline{1,n}$, we get \cite{3}

\begin{proposition} {\normalfont{(a)}}    With respect to the transformation of coordinates $\bar{x}^{i} =\bar{x}^{i} (x^{1} ,...,x^{n} )$, $i=\overline{1,n}$, $\det \left(\frac{\partial \bar{x}^{i} }{\partial x^{j} } \right)\ne 0$, corresponding to the charts $U$, $U'$, $U\cap U'\ne \emptyset $, we have the relations

\begin{equation} \label{6)} d(\bar{x}^{i} )^{\alpha } =\mathop{J_{j}^{i} }\limits^{\alpha } (\bar{x},x)d(x^{j} )^{\alpha } ,\; \; D_{\bar{x}^{i} }^{\alpha } =\mathop{J_{i}^{j} }\limits^{\alpha } (x,\bar{x})D_{x^{j} }^{\alpha } , \end{equation}

where $\mathop{J_{j}^{i} }\limits^{\alpha } (x,\bar{x})=\frac{1}{\Gamma (1+\alpha )} D_{\bar{x}^{j} }^{\alpha } (x^{j} )^{\alpha } .$

{\normalfont{(b)}}    The transformation of coordinates on $(\pi ^{\alpha k} )^{-1} (U\cap U')\subset Osc^{\alpha k} (M)$ are given by
\begin{multline} \label{7)} \frac{\Gamma (\alpha (a-1))}{\Gamma (\alpha )} \bar{y}^{i(\alpha a)} =\Gamma (1+\alpha )\mathop{J_{j}^{i} }\limits^{\alpha } (\bar{y}^{\alpha (a-1)} ,x)y^{j(\alpha )} \\ +\frac{\Gamma (2\alpha )}{\Gamma (\alpha )}  \sum _{b=1}^{a-1}\mathop{J_{j}^{i} }\limits^{\alpha } (\bar{y}^{\alpha (a-1)} ,y^{\alpha b} )y^{j((b+1)\alpha )}    +\frac{\Gamma (\alpha (a-1))}{\Gamma (\alpha )} y^{i(\alpha k)}, \end{multline}

where $a=\overline{2,k}$ and $(x^{i} ,y^{i(\alpha )} ,...,y^{i(\alpha k)} )\in (\pi ^{\alpha k} )^{-1} (U)$.
\end{proposition}

\section{The fractional jet bundle on a differentiable manifold. Geometrical objects}
\subsection{The fractional jet bundle}

 \textit{The fractional jet bundle of order} $\alpha $ on the manifold $M$ is the space $J^{\alpha } ({\mathbb R},\; M)={\mathbb R}\times Osc^{\alpha } (M)$, where $Osc^{\alpha } (M)$ is the fractional osculator bundle. The triplet $(J^{\alpha } ({\mathbb R},\; M),\; \pi _{0}^{\alpha } ,\; M)$ has a structure of differentiable fiber bundle, where $\pi _{0}^{\alpha } :J^{\alpha } ({\mathbb R},\; M)\to M$ is the canonical projection. If $(x^{i} )$, $i=\overline{1,n}$, are the coordinate functions on the chart $U\subset M$, then the coordinate functions on $(\pi _{0}^{\alpha } )^{-1} (U)$ are given by $(t,x^{i} ,y^{i(\alpha )} )$ where $y{}^{i(\alpha )} =\frac{1}{\Gamma (1+\alpha )} D_{t}^{\alpha } x^{i} (t)\left|{}_{t=0} \right. $, $i=\overline{1,n}$.

 From the properties of the Caputo fractional derivative (subsection 2.1) and from Proposition 1, it follows:

\begin{proposition} Let us consider the functions $(t)^{\alpha } $, $(x^{i} )^{\alpha } $, $(y^{i(\alpha )} )^{\alpha } $ of \\ ${\mathcal{F}}((\pi _{0}^{\alpha } )^{-1} (U))$, the 1-forms $\frac{1}{\Gamma (1+\alpha )} d(t)^{\alpha } $, $\frac{1}{\Gamma (1+\alpha )} d(x^{i} )^{\alpha } $, $\frac{1}{\Gamma (1+\alpha )} d((y^{i(\alpha )} )^{\alpha } )$ of \\ ${\mathcal{D}}^{1} ((\pi _{0}^{\alpha } )^{-1} (U))$ and the operators $D_{t}^{\alpha } $, $D_{x^{i} }^{\alpha } $, $D_{y^{i(\alpha )} }^{\alpha } $, $i=\overline{1,n}$. The following relations hold:
\[\begin{array}{l} {D_{t}^{\alpha } \left(\frac{1}{\Gamma (1+\alpha )} t^{\alpha } \right)=1,\; \; D_{x^{i} }^{\alpha } \left(\frac{1}{\Gamma (1+\alpha )} (x^{j} )^{\alpha } \right)=\delta _{i}^{j} ,\; \; D_{y^{i(\alpha )} }^{\alpha } \left(\frac{1}{\Gamma (1+\alpha )} y^{j(\alpha )} \right)=\delta _{i}^{j} ,} \\ {\frac{1}{\Gamma (1+\alpha )} d(t^{\alpha } )(D_{t}^{\alpha } )=1,\; \; \frac{1}{\Gamma (1+\alpha )} d(x^{i} )^{\alpha } (D_{x^{j} }^{\alpha } )=\delta _{j}^{i} ,\;}\\ {\frac{1}{\Gamma (1+\alpha )} d(y^{i(\alpha )} )^{\alpha } (D_{y^{j(\alpha )} }^{\alpha } )=\delta _{j}^{i} .} \end{array}\]
\end{proposition}
 The module generated by the operators $D_{t}^{\alpha } $, $D_{x^{i} }^{\alpha } $, $D_{y^{i(\alpha )} }^{\alpha } $, $i=\overline{1,n}$, will be denoted by ${\rm {\mathfrak X}}^{\alpha } ((\pi _{0}^{\alpha } )^{-1} (U))$. For $\alpha \to 1$ this module represents the module of the vector fields defined on $\pi _{0}^{-1} (U)$.

 Let us consider two charts $U$, $U'$ on $M$ with $U\cap U'\ne \emptyset $, $(\pi _{0}^{\alpha } )^{-1} (U)$, $(\pi _{0}^{\alpha } )^{-1} (U')\subset J^{\alpha } ({\mathbb R},\; M)$ the corresponding charts on $J^{\alpha } ({\mathbb R},\; M)$ and the coordinate functions $(x^{i} )$, $(\bar{x}^{i} )$, respectively, $(t,x^{i} ,y^{i(\alpha )} )$, $(t,\bar{x}^{i} ,\bar{y}^{i(\alpha )} )$. From Proposition 1, we obtain the transformations of coordinates

\begin{equation} \label{ZEqnNum230610} \begin{array}{l} {\bar{x}^{i} =\bar{x}^{i} (x^{1} ,...,x^{n} )} \\ {\bar{y}^{i(\alpha )} =\mathop{J_{j}^{i} }\limits^{\alpha } (x,\bar{x})y^{j(\alpha )} .} \end{array} \end{equation}

\subsection{Geometrical objects on $J^{\alpha } ({\mathbb R},\; M)$ }

 On the manifold $J^{\alpha } ({\mathbb R},\; M)$, the following canonical structures may be defined:

\begin{equation} \label{ZEqnNum749847} \begin{array}{l} {\mathop{\theta }\limits^{\alpha } _{1} =d(t^{\alpha } )\otimes (D_{t}^{\alpha } +y^{i(\alpha )} D_{x^{i} }^{\alpha } )} \\ {\mathop{\theta }\limits^{\alpha } _{2} =\mathop{\theta ^{i} }\limits^{\alpha } \otimes D_{x^{i} }^{\alpha } ,\; \; \mathop{\theta ^{i} }\limits^{\alpha } =\frac{1}{\Gamma (1+\alpha )} (d(x^{i} )^{\alpha } -y^{i(\alpha )} d(t^{\alpha } ))} \\ {\mathop{S}\limits^{\alpha } =\mathop{\theta ^{i} }\limits^{\alpha } \otimes D_{y^{i(\alpha )} }^{\alpha } } \\ {\mathop{V_{i} }\limits^{\alpha } =D_{y^{i(\alpha )} }^{\alpha } .} \end{array} \end{equation}

From \eqref{ZEqnNum230610} and from Proposition 1, it follows that the structures \eqref{ZEqnNum749847} have geometrical character.

 The vector field $\mathop{\Gamma }\limits^{\alpha } \in {\rm {\mathfrak X}}^{\alpha } ((\pi _{0}^{\alpha } )^{-1} (U))$ is called a \textit{fractional vector field} (\textit{FVF}) iff

\begin{equation} \label{10)} d(t)^{\alpha } (\mathop{\Gamma }\limits^{\alpha } )=1,\; \mathop{\; \theta ^{i} }\limits^{\alpha } (\mathop{\Gamma }\limits^{\alpha } )=0,\; \; i=\overline{1,n}. \end{equation}

In local coordinates, (\textit{FVF}) is given by

\begin{equation} \label{11)} \mathop{\Gamma }\limits^{\alpha } =D_{t}^{\alpha } +y^{i(\alpha )} D_{x^{i} }^{\alpha } +F^{i} D_{y^{i(\alpha )} }^{\alpha } , \end{equation}

where $F^{i} \in C^{\infty } ((\pi _{0}^{\alpha } )^{-1} (U))$, $i=\overline{1,n}$. The integral curves of (\textit{FVF}) are the solutions of the fractional differential equations (\textit{FDE})

\begin{equation} \label{ZEqnNum663740} D_{t}^{2\alpha } x^{i} (t)=F^{i} (t,x(t),D_{t}^{\alpha } x(t)),\; \; i=\overline{1,n}. \end{equation}

The system \eqref{ZEqnNum663740}, with given initial conditions, admits solution \cite{5}.

 Let $L\in C^{\infty } (J^{\alpha } ({\mathbb R},\; M))$ be a \textit{fractional Lagrange function}. By definition, the \textit{Cartan} \textit{fractional }1-\textit{form} is the 1-form $\mathop{\theta _{L} }\limits^{\alpha } $ given by

\begin{equation} \label{13)} \mathop{\theta _{L} }\limits^{\alpha } =Ld(t)^{\alpha } +\mathop{S}\limits^{\alpha } (L). \end{equation}

The \textit{Cartan }2-\textit{form }$\mathop{\omega _{L} }\limits^{\alpha } $ is defined by

\begin{equation} \label{14)} \mathop{\omega _{L} }\limits^{\alpha } =d^{\alpha } \mathop{\theta _{L} }\limits^{\alpha } , \end{equation}

where $d^{\alpha } $ is the fractional exterior differential

\begin{equation} \label{15)} d^{\alpha } =d(t)^{\alpha } D_{t}^{\alpha } +d(x^{i} )^{\alpha } D_{x^{i} }^{\alpha } +d(y^{i(\alpha )} )^{\alpha } D_{y^{i(\alpha )} }^{\alpha } . \end{equation}

 In the chart $(\pi _{0}^{\alpha } )^{-1} (U)$, $\mathop{\theta _{L} }\limits^{\alpha } $ and $\mathop{\omega _{L} }\limits^{\alpha } $ are given by

\begin{equation} \label{16)} \begin{array}{l} {\mathop{\theta _{L} }\limits^{\alpha } =(L-\frac{1}{\Gamma (1+\alpha )} y^{i(\alpha )} D_{y^{i(\alpha )} }^{\alpha } (L))d(t)^{\alpha } +\frac{1}{\Gamma (1+\alpha )} D_{y^{i(\alpha )} }^{\alpha } (L)d(x^{i} )^{\alpha } } \\ {\mathop{\omega _{L} }\limits^{\alpha } =A_{i} d(t)^{\alpha } \wedge d(x^{i} )^{\alpha } +B_{i} d(t)^{\alpha } \wedge d(y^{i(\alpha )} )^{\alpha } +} \\ {A_{ij} d(x^{i} )^{\alpha } \wedge d(x^{j} )^{\alpha } +B_{ij} d(x^{i} )^{\alpha } \wedge d(y^{j(\alpha )} )^{\alpha } ,} \end{array} \end{equation}

where

\[\begin{array}{l} {A_{i} =\frac{1}{\Gamma (1+\alpha )} D_{t}^{\alpha } D_{y^{i(\alpha )} }^{\alpha } (L)+\frac{1}{\Gamma (1+\alpha )} y^{j(\alpha )} D_{x^{i} }^{\alpha } D_{y^{j(\alpha )} }^{\alpha } (L)-D_{x^{i} }^{\alpha } (L)} \\ {B_{i} =\frac{1}{\Gamma (1+\alpha )} D_{y^{i(\alpha )} }^{\alpha } (y^{j(\alpha )} D_{y^{j(\alpha )} }^{\alpha } (L))} \\ {A_{ij} =D_{x^{i} }^{\alpha } D_{y^{i(\alpha )} }^{\alpha } (L),\; \; B_{ij} =-D_{y^{j(\alpha )} }^{\alpha } D_{y^{i(\alpha )} }^{\alpha } (L).} \end{array}\]

\begin{proposition} If the fractional Lagrange function is regular i.e., \\ $\det \left(\mathop{g_{ij} }\limits^{\alpha } \right)\ne 0$, $\mathop{g_{ij} }\limits^{\alpha } =D_{y^{i(\alpha )} }^{\alpha } D_{y^{j(\alpha )} }^{\alpha } L$, then there is a fractional vector field (\textit{FVF}) $\mathop{\Gamma _{L} }\limits^{\alpha } $ such that $i_{\mathop{\Gamma _{L} }\limits^{\alpha } } (\mathop{\omega _{L} }\limits^{\alpha } )=0$. In the chart $(\pi _{0}^{\alpha } )^{-1} (U)$, $\mathop{\Gamma _{L} }\limits^{\alpha } $ is given by

\begin{equation} \label{ZEqnNum716035} \mathop{\Gamma _{L} }\limits^{\alpha } =D_{t}^{\alpha } +y^{i(\alpha )} D_{x^{i} }^{\alpha } +\mathop{M^{i} }\limits^{\alpha } D_{y^{i(\alpha )} }^{\alpha } , \end{equation}

where

\begin{equation} \label{18)} \begin{array}{l} {\mathop{M^{i} }\limits^{\alpha } =\mathop{g^{ik} }\limits^{\alpha } (D_{x^{k} }^{\alpha } L-d_{t}^{\alpha } (D_{y^{k(\alpha )} }^{\alpha } L))} \\ {\left(\mathop{g^{ik} }\limits^{\alpha } \right)=\left(\mathop{g_{ik} }\limits^{\alpha } \right)^{-1} ,\; \; d_{t}^{\alpha } =D_{t}^{\alpha } +y^{i(\alpha )} D_{x^{i} }^{\alpha } .} \end{array} \end{equation}
\end{proposition}

\subsection{The fractional Euler-Lagrange equations }

 Let $c:t\in [0,1]\to (x^{i} (t))\in M$ be a parameterized curve such that $Imc\subset U\subset M$. The extension of the curve $c$ to $J^{\alpha } ({\mathbb R},\; M)$ is the curve $c^{\alpha } :t\in [0,1]\to (t,x^{i} (t),y^{i(\alpha )} (t))\in J^{\alpha } ({\mathbb R},\; M)$ with $Imc^{\alpha } \subset (\pi _{0}^{\alpha } )^{-1} (U)\subset J^{\alpha } ({\mathbb R},\; M)$. Let $L\in C^{\infty } (J^{\alpha } ({\mathbb R},\; M))$ be a fractional Lagrange function. The action of $L$ along the curve $c^{\alpha } $ is

\begin{equation} \label{19)} {\it A}(c^{\alpha } )=\int _{0}^{1}L(t,x(t),y^{(\alpha )} (t))dt . \end{equation}

Let $c_{\varepsilon } :t\in [0,1]\to (x^{i} (t,\varepsilon ))\in M$ be a family of curves, with $\varepsilon $ sufficiently small in absolute value so that $Imc_{\varepsilon } \subset M$, $c_{0} (t)=c(t)$, $D_{\varepsilon }^{\alpha } c_{\varepsilon } (0)=D_{\varepsilon }^{\alpha } c_{\varepsilon }(1)=0$. The action of $L$ on the curves $c_{\varepsilon }^{\alpha } $ is

\begin{equation} \label{ZEqnNum399731} {\it A}(c_{\varepsilon }^{\alpha } )=\int _{0}^{1}L(t,x(t,\varepsilon ),y^{(\alpha )} (t,\varepsilon ))dt , \end{equation}

where $y^{i(\alpha )} (t,\varepsilon )=\frac{1}{\Gamma (1+\alpha )} D_{t}^{\alpha } x^{i} (t,\varepsilon )$, $i=\overline{1,n}$. The action \eqref{ZEqnNum399731} has a fractional extremal value if

\begin{equation} \label{21)} D_{\varepsilon }^{\alpha } {\it A}(c_{\varepsilon }^{\alpha } )\left|_{\varepsilon =0} \right. =0. \end{equation}

Using the properties of the Caputo fractional derivative (subsection 2.1), it results

\begin{proposition} {\normalfont{(a)}} If the action \eqref{ZEqnNum399731} reaches a fractional extremal value then a necessary condition is that $c(t)$ satisfies \textit{the fractional Euler-Lagrange equations}

\begin{equation} \label{ZEqnNum706067} \begin{array}{l} {D_{x^{i} }^{\alpha } L-d_{t}^{2\alpha } (D_{y^{i(\alpha )} }^{\alpha } L)=0,\; \; i=\overline{1,n}} \\ {d_{t}^{2\alpha } =D_{t}^{\alpha } +y^{i(\alpha )} D_{x^{i} }^{\alpha } +y^{i(2\alpha )} D_{y^{i(\alpha )} }^{\alpha } .} \end{array} \end{equation}

{\normalfont{(b)}} If the fractional Lagrange function is nondegenerated, then the equations \eqref{ZEqnNum706067} are the fractional differential equations associated to the fractional vector field $\mathop{\Gamma _{L} }\limits^{\alpha } $ given by \eqref{ZEqnNum716035}.

{\normalfont{(c)}} If the fractional Lagrange function is nondegenerated, then the system \eqref{ZEqnNum706067} may be written in the form of \textit{the fractional Hamilton equations}

\begin{equation} \label{23)} D_{t}^{\alpha } p_{i}^{(\alpha )} =-D_{x^{i} }^{\alpha } H,\; \; D_{t}^{\alpha } x^{i} =D_{p_{i}^{(\alpha )} }^{\alpha } H, \end{equation}

where

\begin{equation} \label{24)} \begin{array}{l} {H=p_{i}^{(\alpha )} D_{t}^{\alpha } x^{i} -L(t,x(t),y^{(\alpha )} (t))} \\ {p_{i}^{(\alpha )} =D_{y^{i(\alpha )} }^{\alpha } L(t,x,y^{(\alpha )} ),\; \; i=\overline{1,n}.} \end{array} \end{equation}

{\normalfont{(d)}} If for $f,h:J^{1} ({\mathbb R},\; M)^{*} \to {\mathbb R}$ \textit{the fractional Poisson bracket} is defined by

\begin{equation} \label{25)} \{ f,h\} ^{\alpha } =D_{p_{i}^{(\alpha )} }^{\alpha } fD_{x^{i} }^{\alpha } g-D_{x^{i} }^{\alpha } fD_{p_{i}^{(\alpha )} }^{\alpha } g, \end{equation}

where the local coordinates on $J^{1} ({\mathbb R},\; M)^{*} $ are $(x,p^{(\alpha )} )$, then

\begin{equation} \label{26)} \{ H,p_{i}^{(\alpha )} \} ^{\alpha } =D_{t}^{\alpha } p_{i}^{(\alpha )} ,\; \; \; \; \{ H,x^{i} \} ^{\alpha } =D_{t}^{\alpha } x^{i} ,\; \; i=\overline{1,n}. \end{equation}

\end{proposition} 
\section{Economic models described by fractional differential equations }
\subsection{The fractional model of Liviatan-Samuelson }

 Let us consider the fractional Lagrange function $L\in {\mathcal{F}}(J^{\alpha } ({\mathbb R},\; M))$ given by

\begin{equation} \label{ZEqnNum470483} L(t,x,y^{(\alpha )} )=L_{1} (x,y^{(\alpha )} )E_{\alpha } (-\rho t^{\alpha } ), \end{equation}

where $E_{\alpha } $ is the Mittag-Leffler function, $E_{\alpha } (t)=\sum _{k=0}^{\infty }\frac{t^{\alpha k} }{\Gamma (1+\alpha k)}  $ and $\rho >0$ is the discount rate. Using the relation $D_{t}^{\alpha } E_{\alpha } (-\rho t^{\alpha } )=-\rho E_{\alpha } (-\rho t^{\alpha } )$ and Proposition 4 we obtain

\begin{proposition}{\normalfont{(a)}} The fractional Euler-Lagrange equations \eqref{ZEqnNum706067} for \eqref{ZEqnNum470483} are
\begin{equation} \label{28)} y^{j(2\alpha )} D_{y^{i(\alpha )} }^{\alpha } D_{y^{j(\alpha )} }^{\alpha } L_{1} +y^{j(\alpha )} D_{x^{j} }^{\alpha } D_{y^{i(\alpha )} }^{\alpha } L_{1} -\rho D_{y^{i(\alpha )} }^{\alpha } L_{1} -D_{x^{i} }^{\alpha } L_{1} =0,\; \; i=\overline{1,n}. \end{equation}

{\normalfont{(b)}} If $L_{1} \in {\mathcal{F}}(J^{\alpha } ({\mathbb R},\; {\mathbb R}))$ is of the form

\begin{equation} \label{29)} L_{1} (x,y^{(\alpha )} )=U(g(x)-y^{(\alpha )} ), \end{equation}

where $U$ is the utility (welfare) function and $c=g(x)-y^{(\alpha )} $ is the consumption function, then the fractional Euler-Lagrange equation is

\begin{equation} \label{30)} \begin{array}{l} {\Gamma (1+\alpha )^{2} U''(g(x)-y^{(\alpha )} )(y^{(\alpha )} )^{2\alpha } -\Gamma (1+\alpha )D_{x}^{\alpha } g(x)U''(g(x)-y^{(\alpha )} )y^{(\alpha )} } \\ {+\rho U'(g(x)-y^{(\alpha )} )\Gamma (1+\alpha )-U'(g(x)-y^{(\alpha )} )D_{x}^{\alpha } g(x)=0.} \end{array} \end{equation}

\end{proposition}

\begin{proposition} If $L_{1} \in {\mathcal{F}}(J^{\alpha } ({\mathbb R},\; {\mathbb R}))$ is given by
\begin{equation} \label{31)} L_{1} (x,y^{(\alpha )} )=-a_{1} (y^{(\alpha )} )^{2\alpha } -a_{2} (y^{(\alpha )} )^{\alpha } x^{\alpha } -a_{3} x^{2\alpha } ,\; \; a_{1} ,a_{2} ,a_{3} \in {\mathbb R}, \end{equation}

then the fractional Euler-Lagrange equation is
\begin{equation} \label{ZEqnNum454786} \begin{array}{l} {a_{1} \Gamma (1+\alpha )\Gamma (1+2\alpha )y^{(2\alpha )} -(a_{2} \Gamma (1+\alpha )^{2} +\rho a_{1} \Gamma (1+2\alpha ))(y^{(\alpha )} )^{\alpha }}\\ {+a_{2} \Gamma (1+\alpha )^{3} y^{(\alpha )}  -(a_{3} \Gamma (1+2\alpha )+\rho a_{2} \Gamma (1+\alpha )^{2} )x^{\alpha } =0.} \end{array} \end{equation}
\end{proposition}

If $a_{1} =a_{3} =\frac{1}{2} $, $a_{2} =a$ and $\alpha \to 1$, then the function $L_{1} $ and the equation \eqref{ZEqnNum454786} become

\begin{equation} \label{ZEqnNum368526} L{}_{1} (x,\dot{x})=-\frac{1}{2} \dot{x}^{2} -ax\dot{x}-\frac{1}{2} x^{2} . \end{equation}

The equation \eqref{ZEqnNum368526} represents the classic model of Samuelson \cite{7}.

\subsection{Fractional economic models with restrictions }

 Let us consider the Lagrange function $L\in C^{\infty } (J^{\alpha } ({\mathbb R},\; M))$ and the function $F\in C^{\infty } (J^{\alpha } ({\mathbb R},\; M))$. The fractional Euler-Lagrange equations of $L$ on the restriction $F(x,y^{(\alpha )} )=0$, $(x,y^{(\alpha )} )\in (\pi _{0}^{\alpha } )^{-1} (U)$, are given by the fractional Euler-Lagrange equations of the fractional Lagrange function

\begin{equation} \label{ZEqnNum394608} L_{2} (t,\lambda ,x,y^{(\alpha )} )=L(t,x,y^{(\alpha )} )+\lambda F(x,y^{(\alpha )} ), \end{equation}

where $\lambda (t)$ is a Lagrange multiplier. From Proposition 4, we obtain

\begin{proposition} {\normalfont{(a)}} The fractional Euler-Lagrange equations of \eqref{ZEqnNum394608} are:

\begin{equation} \label{ZEqnNum220801} \begin{array}{l} {D_{x^{i} }^{\alpha } L+\lambda D_{x^{i} }^{\alpha } F-d_{t}^{2\alpha } (D_{y^{i(\alpha )} }^{\alpha } L)-\lambda y^{j(\alpha )} D_{x^{j} }^{\alpha } (D_{y^{i(\alpha )} }^{\alpha } F)} \\ {-\lambda y^{j(\alpha )} D_{y^{j(\alpha )} }^{\alpha } (D_{y^{i(\alpha )} }^{\alpha } F)-D_{t}^{\alpha } \lambda D_{y^{i(\alpha )} }^{\alpha } F=0,\; \; i=\overline{1,n.}} \end{array} \end{equation}

{\normalfont{(b)}} If the Lagrange function is given by \eqref{ZEqnNum470483} then the fractional Euler-Lagrange equations \eqref{ZEqnNum220801} become

\begin{equation} \label{ZEqnNum221340}
\begin{array}{l} {E_{\alpha } (-\rho t^{\alpha } )(D_{x^{i} }^{\alpha } L_{1} +\rho D_{y^{i(\alpha )} }^{\alpha } L_{1} -y^{j(\alpha )} D_{x^{j} }^{\alpha } (D_{y^{i(\alpha )} }^{\alpha } L_{1} )} \\ {-y^{j(2\alpha )} D_{y^{j(\alpha )} }^{\alpha } (D_{y^{i(\alpha )} }^{\alpha } L_{1} ))+\lambda (D_{x^{i} }^{\alpha } F-y^{j(\alpha )} D_{x^{j} }^{\alpha } (D_{y^{i(\alpha )} }^{\alpha } F)} \\ {-y^{j(2\alpha )} D_{y^{j(\alpha )} }^{\alpha } (D_{y^{i(\alpha )} }^{\alpha } F))-D_{t}^{\alpha } \lambda D_{y^{i(\alpha )} }^{\alpha } F=0,} \end{array}
 \end{equation}

for $i=\overline{1,n}$.
\end{proposition} 
 The fractional model of investments with restriction is described by the function $L_{1} (K,I,N)$ where $K(t)=x^{1} (t)$, $I(t)=x^{2} (t)$, $N(t)=x^{3} (t)$ represent the capital stock, the investment and the labor, respectively. The restriction is given by $F(K^{(\alpha )} ,K,I,N)=\phi (K,I,N)-K^{(\alpha )} =0,$ $K^{(\alpha )} =D_{t}^{\alpha } K$. From \eqref{ZEqnNum221340} we obtain the fractional Euler-Lagrange equations

\begin{equation} \label{ZEqnNum417497} \begin{array}{l} {E_{\alpha } (-\rho t^{\alpha } )D_{K}^{\alpha } L_{1} +\lambda D_{K}^{\alpha } \phi =-D_{t}^{\alpha } \lambda } \\ {E_{\alpha } (-\rho t^{\alpha } )D_{I}^{\alpha } L_{1} +\lambda D_{K}^{\alpha } \phi =0} \\ {E_{\alpha } (-\rho t^{\alpha } )D_{N}^{\alpha } L_{1} +\lambda D_{N}^{\alpha } \phi =0.} \end{array} \end{equation}

If $L_{1} $ and $\phi $ satisfy the relations

\begin{equation} \label{38)} \begin{array}{l} {K^{\alpha } D_{K}^{\alpha } L_{1} +I^{\alpha } D_{I}^{\alpha } L_{1} +N^{\alpha } D_{N}^{\alpha } L{}_{1} =\frac{1}{\Gamma (1+\alpha )} rL_{1} ,\; \; r\in {\mathbb R}} \\ {K^{\alpha } D_{K}^{\alpha } \phi +I^{\alpha } D_{I}^{\alpha } \phi +N^{\alpha } D_{N}^{\alpha } \phi =\frac{1}{\Gamma (1+\alpha )} \phi ,} \end{array} \end{equation}

from \eqref{ZEqnNum417497} we get

\begin{equation} \label{39)} rE_{\alpha } (-\rho t^{\alpha } )L_{1} =-\lambda K^{(\alpha )} -\Gamma (1+\alpha )K^{\alpha } \lambda ^{(\alpha )} . \end{equation}

For $\alpha \to 1$ the classic model of investments \cite{7} is obtained.

\end{document}